\documentclass{article}
\usepackage[utf8]{inputenc}

\usepackage{subfigure}
\usepackage{tikz}

\usepackage[utf8]{inputenc}

\usepackage[utf8]{inputenc}
\usepackage[english]{babel}
\usepackage{amsmath}
\usepackage{times}
\usepackage{graphicx}
\usepackage{float}
\usepackage[margin=1in]{geometry}
\usepackage{blindtext}
\usepackage{amssymb}
\usepackage{tikz}
\usepackage{amsthm}
\usepackage{mathtools}
\usepackage{textcmds}

\usepackage{comment}
\newtheorem{theorem}{Theorem}[section]
\newtheorem{corollary}[theorem]{Corollary}

\newtheorem{lemma}[theorem]{Lemma}
\newtheorem{conjecture}[theorem]{Conjecture}
\theoremstyle{remark}

\theoremstyle{definition}

\date{}
\begin{document}

\title{Spectral extremal graphs for intersecting cliques}

\author{
Dheer Noal Desai\thanks{Department of Mathematical Sciences, University of Delaware, U.S.A. \texttt{dheernsd@udel.edu}} \and
Liying Kang\thanks{Department of Mathematics, Shanghai University,
Shanghai 200444, P.R. China \texttt{lykang@shu.edu.cn}. This work is partially supported by the National
Nature Science Foundation of China (grant numbers 11871329, 11971298).} \and
Yongtao Li\thanks{School of Mathematics, Hunan University, Changsha, P.R. China \texttt{ytli0921@hnu.edu.cn}. This work is partially supported by NSFC (Grant Nos. 11931002, 11671124).} \and
Zhenyu Ni\thanks{Department of Mathematics, Shanghai University,
Shanghai 200444, P.R. China \texttt{1051466287@qq.com}}\and
Michael Tait\thanks{Department of Mathematics and Statistics, Villanova University, U.S.A. \texttt{michael.tait@villanova.edu}. Research is partially supported by National Science Foundation grant DMS-2011553.} \and
Jing Wang\thanks{Department of Mathematics, Shanghai University,
Shanghai 200444, P.R. China \texttt{wj517062214@163.com}} \\[2ex]
}

\maketitle

\vspace{-0.5cm}

\begin{abstract}
The $(k,r)$-fan is the graph consisting of $k$ copies of the
complete graph $K_r$ which intersect in a single vertex, and is denoted by $F_{k,r}$.
Erd\H{o}s, F\"uredi, Gould and Gunderson [J. Combin. Theory Ser. B 64 (1995) 89--100]
determined the maximum number of edges in an $n$-vertex graph that does not contain $F_{k,3}$ as a subgraph. Furthermore, Chen, Gould, Pfender and Wei
[J. Combin. Theory Ser. B 89 (2003) 159--171] proved the analogous result on $F_{k,r}$ for the general case  $r\ge 3$.
In this paper, we show that for sufficiently large $n$, the
graphs of order $n$ that contain no copy of $F_{k,r}$ and
attain the maximum spectral radius
 are also edge-extremal. That is, such graphs must have $\mathrm{ex}(n, F_{k,r})$ edges.
 \end{abstract}


\section{Introduction}
In this paper, we consider only simple and undirected graphs. Let $G$ be a simple connected  graph with vertex set $V(G)=\{v_1, \ldots, v_n\}$ and edge set $E(G)=\{e_1, \ldots, e_m\}$.
For a vertex $v\in V(G)$, we write $N(v)$  for the set of neighbors of $v$.
Let $d(v)$  be the degree of a vertex $v$ in $G$. That is, $d(v)=|N(v)|$.
Let $S$ be a set of vertices. We write
$N_S(v)$ for the set of neighbors of $v$ in the set $S$, and
$d_S(v)$
for the number of neighbors of $v$ in the set $S$,
that is, $d_S(v)=|N_S(v)|=|N(v)\cap S|$.
And we denote by $e(S)$ the number of edges contained in $S$.

The main tasks in extremal graph theory
are to maximize or minimize a graph parameter over a specific family of graphs.
 The {\em Tur\'an number} of a graph $F$ is the maximum number of edges that may be
  in an $n$-vertex graph without a subgraph isomorphic to $F$, and
 this quantity is usually  denoted by $\mathrm{ex}(n, F)$.
  We say that a graph $G$ is $F$-free if it does not contain a subgraph
 isomorphic to $F$, i.e., $G$ contains no copy of $F$.
  A graph on $n$ vertices with no subgraph $F$ and with $\mathrm{ex}(n, F)$ edges is called an {\em extremal graph} for $F$ and we denote by $\mathrm{Ex}(n, F)$ the set of all extremal graphs on $n$ vertices for $F$.
  It is  a cornerstone of extremal graph theory
  to
  investigate both $\mathrm{ex}(n, F)$ and $\mathrm{Ex}(n, F)$ for various graphs $F$;
  see \cite{FS13, Keevash11,Sim13} for related surveys.

Dating back to 1941,
Tur\'{a}n \cite{Turan41} first raised the  natural question of determining
 $\mathrm{ex}(n,K_{r+1})$ where $K_{r+1}$ is the complete graph on
  $r+1$ vertices.
 Let $T_r(n)$ denote the complete $r$-partite graph on $n$ vertices where
 its part sizes are as equal as possible, i.e., each part has size $\lfloor n/r \rfloor $ or $\lceil n/r \rceil$.
Tur\'{a}n \cite{Turan41}
  extended a result of Mantel \cite{Man07}
and  obtained that if $G$ is an $n$-vertex graph containing no $K_{r+1}$,
then $e(G)\le e(T_r(n))$, equality holds if and only if $G=T_r(n)$.
There are many extensions and generalizations of Tur\'{a}n's result;
see, e.g., \cite[p. 294]{Bollobas78}.
After this result, the problem of determining $\mathrm{ex}(n, F)$ is usually called the
Tur\'{a}n-type extremal problem.
 The most celebrated result is a theorem of
 Erd\H{o}s, Stone and Simonovits \cite{ES46,ES66}, which  states that
 \begin{equation} \label{eqESS}
  \mathrm{ex}(n,F) = \left( 1- \frac{1}{\chi (H) -1}  \right)
  \frac{n^2}{2} + o(n^2),
  \end{equation}
 where $\chi (F)$ is the vertex-chromatic number of $H$.
 This provides good asymptotic estimates for the extremal numbers of non-bipartite graphs.
 However, for bipartite graphs, where $\chi (F)=2$, it only gives the bound
 $\mathrm{ex}(n,F)=o(n^2)$.
 The history of studying bipartite graphs began in 1954 with
 the K\H{o}v\'ari--S\'{o}s--Tur\'{a}n theorem \cite{KST54}, which asserts
 that if $K_{s,t}$ is the complete bipartite graph with vertex classes of size $s\ge t$,
 then $\mathrm{ex}(n,K_{s,t})=O(n^{2-1/t})$; see \cite{Furedi96,Furedi96b}
 for more details.
  Although there have been numerous attempts to find better bounds
 of $\mathrm{ex}(n,F)$ for various bipartite graphs $F$, we know very little in this case.
We refer the interested reader to
  the comprehensive survey by F\"{u}redi and Simonovits \cite{FS13}.

\subsection{Background and motivation}
In this section, we shall review
 the exact value of $\mathrm{ex}(n,F)$ for some
special graphs $F$, instead of the asymptotic estimation.
A graph on $2k+1$ vertices consisting of $k$ triangles which intersect in exactly one
common vertex is called a {\it $k$-fan} (also known as the friendship graph)
and denoted by $F_k$.
Since $\chi (F_k)=3$, the Erd\H{o}s--Stone--Simonovits
theorem in (\ref{eqESS}) implies that
$\mathrm{ex}(n,F_k)= n^2/4 + o(n^2)$.
In 1995, Erd\H{o}s, F\"uredi, Gould and Gunderson \cite{Erdos95} proved  the following exact result.

\begin{theorem}[Erd\H{o}s et al. \cite{Erdos95}] \label{thmErdos95} 
For every $k \geq 1$, and for every $n\geq 50k^2$, we have
\[\mathrm{ex}(n, F_k)= \left\lfloor \frac {n^2}{4}\right \rfloor+ \left\{
  \begin{array}{ll}
   k^2-k, \quad~~  \mbox{if $k$ is odd,} \\
    k^2-\frac32 k, \quad \mbox{if $k$ is even}.
  \end{array}
\right. \]
\end{theorem}

 \medskip

 A graph on $(r-1)k+1$ vertices consisting of
 $k$ cliques each with $r$ vertices,
 which intersect in exactly one common vertex,
 is called a  $(k,r)$-fan and denoted by $F_{k,r}$.
 Clearly, when $r=3$, $F_{k,3}$ reduces to the general $k$-fan graph
 $F_k$.
 Note that $\chi (F_{k,r})= r$.
Similarly,  the Erd\H{o}s--Stone--Simonovits theorem
also implies that
 $\mathrm{ex}(n,F_{k,r})=(1- \frac{1}{r-1})\frac{n^2}{2} + o(n^2)
 =t_{r-1}(n) + o(n^2)$.
 In 2003, Chen, Gould, Pfender and Wei \cite{Chen03}
proved an exact answer and generalized
 Theorem \ref{thmErdos95} as follows.

 \begin{theorem}[Chen et al. \cite{Chen03}] \label{thmChen}
 For every $k\ge 1$ and $r\ge 2$, if $n\ge 16k^3r^8$, then
 \[ \mathrm{ex}(n,F_{k,r}) = t_{r-1}(n) +
 \left\{
  \begin{array}{ll}
   k^2-k, \quad~~  \mbox{if $k$ is odd,} \\
    k^2-\frac32 k, \quad \mbox{if $k$ is even}.
  \end{array}
\right.   \]
 \end{theorem}

 The extremal graphs of Theorem \ref{thmChen}, denoted by $G_{n,k,r}$, are constructed by taking the $(r-1)$-partite Tur\'{a}n graph $T_{r-1}(n)$ and embedding a graph $G_0$ in one vertex part. If $k$ is odd,  $G_{0}$ is isomorphic to two vertex disjoint copies of $K_k$. If $k$ is even,  $G_{0}$ may be isomorphic to any graph with $2k-1$ vertices, $k^2-\frac{3}{2}k$ edges with maximum degree $k-1$.

\subsection{Spectral extremal problem}

The \emph{adjacency matrix} of $G$ is defined as
$A(G)=(a_{ij})_{n \times n}$ with $a_{ij}=1$ if two vertices $v_i$ and $v_j$ are adjacent in $G$, and $a_{ij}=0$ otherwise.
The spectral radius of $A(G)$ is defined as the largest value among the absolute values of eigenvalues of $A(G)$.
Note that the spectral radius is not necessarily
an eigenvalue.
The celebrated Perron--Frobenius theorem
implies that the spectral radius of $A(G)$ is
a largest eigenvalue  since $A(G)$ is a nonnegative matrix.
The eigenvalues of a graph $G$ are defined as
the eigenvalues of adjacency matrix $A(G)$.
We write $\lambda (G)$ or $\lambda_1(G)$ for the spectral
radius of $G$.
The spectral radius of a graph may at times give some information
about the structure of graphs.
For example, it is well-known  that $\lambda (G)$ is
located between
the average degree  and the maximum degree of $G$,
and the vertex-chromatic number is at most $\lambda (G)+1$;
see \cite[p. 34]{Bapat14} for more details.

In this paper we consider spectral analogues
of Tur\'{a}n-type problems for graphs,
that is, determining the maximum value of eigenvalues
instead of the number of edges among all $n$-vertex $F$-free graphs.
We denote
\[  \mathrm{ex}_{sp}(n,F)= \max \{\lambda (G) :
|G|=n, \text{$G$ is $F$-free} \}. \]
These problems are commonly based on the techniques
applying the eigenvalues or eigenvectors of a graph.
The  fundamental inequality $2e(G)/n \le \lambda (G)$ yields the
following relation:
\begin{equation*}
 \mathrm{ex}(n,F) \le \frac{n}{2} \mathrm{ex}_{sp}(n,F) .
 \end{equation*}
 The problems of studying $\mathrm{ex}_{sp}(n,F)$
 has a rapid development in spectral extremal graph theory recently.
For most graphs, this study is again fairly complete
due in large part to a longstanding work of Nikiforov \cite{NikifSurvey}.
For example, he extended the classical theorem of Tur\'{a}n,
by determining the maximum spectral radius of
any $K_{r+1}$-free graph $G$ on $n$ vertices.

The following problem regarding the adjacency  spectral radius
was proposed in \cite{NikiforovTuran}:
What is the maximum spectral radius of a graph $G$ on $n$
vertices without a subgraph isomorphic to a given graph $F$?
Wilf \cite{Wilf86} and Nikiforov \cite{NikiforovTuran}
obtained spectral strengthening of Tur\'an's theorem
when the forbidden  substructure is the complete graph.
Soon after,
Nikiforov \cite{Nikiforov07} showed that if $G$ is a $K_{r+1}$-free graph on $n$ vertices,
then $\lambda (G)\le \lambda (T_r(n))$, with equality if and only if
$G=T_r(n)$.
Moreover,
Nikiforov \cite{Nikiforov07}
and Zhai and Wang \cite{ZW12}
determined the maximum spectral radius
of $K_{2,2}$-free graphs.
Furthermore, Nikiforov \cite{NikiforovKST}, Babai and Guiduli \cite{BG09}
independently  obtained  the spectral generalization of the K\H{o}vari-S\'os-Tur\'an theorem when the forbidden graph is the complete bipartite graph $K_{s,t}$.
Finally, Nikiforov \cite{NikiforovLAA10} characterized
the spectral radius of graphs without paths and cycles of specified length.
In addition, Fiedler and Nikiforov  \cite{FiedlerNikif} obtained tight sufficient conditions for
graphs to be Hamiltonian or traceable.
For many other spectral analogues of results in extremal graph theory
we refer the reader to the survey  \cite{NikifSurvey}.
It is worth mentioning that a corresponding spectral strengthening  \cite{Niki09b} of
the  Erd\H{o}s--Stone--Simonovits theorem states that
\[  \mathrm{ex}_{sp}(n,F) = \left( 1- \frac{1}{\chi (F)-1}  \right)n + o(n). \]
From this result, we know that
$\mathrm{ex}_{sp}(n,F_k) = n/2 +o(n)$ where $F_k$ is the $k$-fan graph.
Recently,  Cioab\u{a}, Feng,  Tait and  Zhang \cite{CFTZ20}   generalized this bound by
improving the error term $o(n)$ to $O(1)$,
and obtained a spectral counterpart of Theorem \ref{thmErdos95}.
More precisely, they showed that the extremal graphs that
attain the maximum spectral radius in a graph on $n$ vertices containing no
copy of $k$-fan must be in $\mathrm{Ex}(n,F_k)$ for $n$ sufficiently large.

 \begin{theorem}[Cioab\u{a} et al. \cite{CFTZ20}]   \label{thmCFTZ20}
Let $G$ be a graph of order $n$ that does not contain a copy of a $k$-fan, $k \geq 2$.
For sufficiently large $n$, if $G$ has the maximal spectral radius, then
$$G \in \mathrm{Ex}(n, F_k).$$
\end{theorem}

Recall that $F_{k,r}$ is the graph consisting of $k$ cliques
 of order $r$ which intersect in
exactly one common vertex.
In this paper, we shall prove the following theorem, which is an extension
of Theorem \ref{thmCFTZ20}.
\begin{theorem}[Main result] \label{main result}
Let $G$ be a graph of order $n$ that does not contain a copy of $F_{k,r}$, where $k\ge 1$ and $r \geq 2$.
For sufficiently large $n$, if $G$ has the maximal spectral radius, then
$$G \in \mathrm{Ex}(n, F_{k,r}).$$
\end{theorem}

Our theorem  is a spectral result of
the Tur\'{a}n extremal problem for $F_{k,r}$, it
 not only can be viewed as
 an extension of Theorem \ref{thmCFTZ20}, but also
a spectral analogue of Theorem \ref{thmChen}.
Our treatment strategy of the proof is mainly
based on the stability method.
To some extent, this paper could be regarded
as a continuation and development of \cite{CFTZ20}. However, we highlight that there are some differences in the approach compared from \cite{CFTZ20}. In \cite{CFTZ20}, the extremal graph is constant edit distance from a bipartite graph. One of the key steps is to show that the extremal graph has a large bipartite subgraph. To do so, the authors prove a lemma (\cite{CFTZ20} Lemma 7) that relates the number of edges to the spectral radius and the number of triangles in the graph and then use the triangle removal lemma and a stability theorem of F\"uredi \cite{Furedi2015}. Unfortunately, for this problem the extremal graph is constant edit distance from an $(r-1)$-partite graph and for $r>3$ the same approach fails. Instead we use a spectral stability theorem of Nikiforov (See Section \ref{sec: lemmas}).

\medskip

{\it Remark.}
The ideas of this paper were developed independently and
 simultaneously by two groups.
 Since the arguments of our two
papers were similar, we present them as a joint work.

\section{Some Lemmas}\label{sec: lemmas}

In this section,
we state some lemmas which are needed in our proof.

\begin{lemma}[Nikiforov \cite{Niki09JGT}]  \label{lemniki}
Let $r\ge 2, 1/\ln n < c < r^{-8(r+21)(r+1)}, 0< \varepsilon < 2^{-36}r^{-24}$
and $G$ be a graph on $n$ vertices. If $\lambda (G) > (1- \frac{1}{r} - \varepsilon )n$, then one of the following statements holds: \\
(a) $G$ contains a $K_{r+1}(\lfloor c\ln n\rfloor , \ldots ,\lfloor c\ln n\rfloor,
\lceil n^{1-\sqrt{c}}\rceil)$; \\
(b) $G$ differs from $T_r(n)$ in fewer than $(\varepsilon^{1/4} +
c^{1/(8r+8)})n^2$ edges.
\end{lemma}

From the above theorem,
one can easily get the following spectral analogue of the classical
Erd\H{o}s-Simonovits stability theorem \cite{Sim66,Furedi2015}.

\begin{corollary}  \label{coro22}
Let $F$ be a graph with chromatic number $\chi (F)=r+1$.
For every $\varepsilon >0$, there exist $\delta >0$
and $n_0$ such that
if  $G$ is an $F$-free graph on $n\ge n_0$ vertices  with
$\lambda (G) \ge (1- \frac{1}{r} -\delta )n$, then
$G$ can be obtained from $T_r(n)$ by adding and deleting at most
$\varepsilon n^2$ edges.
\end{corollary}

Let $G$ be a simple graph with matching number $\beta(G)$ and maximum degree $\Delta(G)$. For given two integers $\beta$ and $\Delta$, define $f(\beta, \Delta)
=\max\{e(G): \beta(G)\leq \beta, \Delta(G)\leq \Delta \}$.

In 1976, Chv\'atal and Hanson \cite{Chvatal76} obtained the following result.

\begin{lemma}[Chv\'atal-Hanson \cite{Chvatal76}]\label{Chvatal76}
For every two integers $\beta \geq 1$ and $\Delta \geq 1$, we have
$$f(\beta, \Delta)= \Delta \beta +\left\lfloor\frac{\Delta}{2}\right\rfloor
 \left \lfloor \frac{\beta}{\lceil{\Delta}/{2}\rceil }\right \rfloor
 \leq \Delta \beta+\beta.$$
\end{lemma}

We will frequently use a special case proved by Abbott, Hanson and Sauer \cite{Abbott72}:
$$f(k-1,k-1) = \left\{
  \begin{array}{ll}
   k^2-k, \quad~~  \mbox{if $k$ is odd,} \\
    k^2-\frac32 k, \quad  \mbox{if $k$ is even}.
  \end{array}
\right.$$
Furthermore, the extremal graphs attaining the equality case
are exactly those we embedded into the Tur\'{a}n graph $T_{r-1}(n)$  to obtain the extremal $F_{k,r}$-free graph.

Denote by  $K_{n_1,n_2,\ldots,n_{r-1}}$  the complete $(r-1)$-partite graph on $n=\sum_{i=1}^{r-1}n_i$ vertices. For convenience, we assume that $n_1\geq n_2\geq \ldots \geq n_{r-1}>0$.
It is well-known \cite[p. 74]{CDS1980} or \cite{Delorme}
that the characteristic polynomial of $K_{n_1,n_2,\ldots,n_{r-1}}$ is given as
$$\phi(K_{n_1,n_2,\ldots,n_{r-1}}, x)=x^{n-r+1}\left(1-\sum_{i=1}^{r-1}\frac{n_i}{x+n_i}\right)\prod_{j=1}^{r-1}(x+n_j).$$
So the spectral radius $\lambda(K_{n_1,n_2,\ldots,n_{r-1}})$ satisfies the following equation:
\begin{equation}
\sum_{i=1}^{r-1}\frac{n_i}{\lambda(K_{n_1,n_2,\ldots,n_{r-1}})+n_{i}}=1 \label{eq0}
\end{equation}

Feng, Li and Zhang \cite[Theorem 2.1]{FLZ2007} proved implicitly
the following lemma, which can also be seen in
Stevanovi\'{c}, Gutnam and Rehman \cite{Stevanovicetal}.

\begin{lemma}[Feng et al. \cite{FLZ2007}, Stevanovi\'{c} et al. \cite{Stevanovicetal}]\label{rpartitegraph}
If $n_i-n_j\geq 2$, then $$\lambda(K_{n_1,\ldots,n_i-1,\ldots,n_j+1,\ldots,n_{r-1}})>\lambda(K_{n_1,\ldots,n_i,\ldots,n_j,\ldots,n_{r-1}}).$$
\end{lemma}

For a connected graph $G$ on $n$ vertices, let $\mathbf{x}=(\mathbf{x}_1,\ldots,\mathbf{x}_n)^{\mathrm{T}}$ be an eigenvector of $A(G)$ corresponding to $\lambda(G)$. By the celebrated Perron--Frobenius theorem, we can
choose $\mathbf{x}$ as a positive real vector.
\begin{equation}\label{eigenequation}
\lambda(G)\mathbf{x}_i= \sum_{j=1}^n a_{ij}\mathbf{x}_j =
\sum_{j\in N_G(i)}\mathbf{x}_j, \text{ for any } i\in [n].
\end{equation}

Another useful result concerns the Rayleigh quotient:
\begin{equation}\label{Rayleigh}
\lambda(G)=\max_{\mathbf{x}\in \mathbb{R}^{n}} \frac{\mathbf{x}^{\mathrm{T}}A(G)\mathbf{x}}{
\mathbf{x}^{\mathrm{T}}\mathbf{x}}=\max_{\mathbf{x}\in \mathbb{R}^{n}}\frac{2\sum_{\{i,j\}\in E(G)}\mathbf{x}_i\mathbf{x}_j}{\mathbf{x}^{\mathrm{T}}\mathbf{x}}.
\end{equation}

Let $G$ be a graph with a partition of the vertices into $r-1$ non-empty parts $V(G)=V_1\cup V_2\cup \ldots \cup V_{r-1}$. Let $E_{cr}(G)=\cup_{1\leq i<j\leq r-1}E(V_i,V_j)$ be the crossing edges of $G$.
The following lemma was proved in Chen et al. \cite{Chen03}.
\begin{lemma}[Chen et al. \cite{Chen03}]\label{chen}
Suppose $G$ is partitioned as above so that the following conditions are satisfied
\begin{eqnarray}
& &\sum_{j\neq i}\beta(G[V_j])\leq k-1\ \ \  \mbox{and} \ \ \ \Delta(G[V_i])\leq k-1,\label{eqn1}\\
& &d_{G[V_i]}(v)+\sum_{j\neq i}\beta(G[N(v)\cap V_j])\leq k-1,\label{eqn2}
\end{eqnarray}
for any $i\in [r-1]$ and  $v\in V_i$.
If $G$ is $F_{k,r}$-free, then $$\sum_{i=1}^{r-1}|E(G[V_i])|-\left(\sum_{1\leq i<j\leq r-1}|V_i||V_j|-|E_{cr}(G)|\right)\leq f(k-1,k-1).$$
\end{lemma}

 \section{Proof of Theorem \ref{main result}}

 In the sequel,
we always assume that $G$ is
 a graph on $n$ vertices
 containing no $F_{k,r}$ as a subgraph and attaining
 the maximum spectral radius.
 The aim of this section is to prove that
 $e(G) = \mathrm{ex}(n, F_{k,r})$ for $n$ large enough.

First of all, we note that $G$ must be connected
since adding an edge between different
components will increase the spectral radius
and also keep $G$ being $F_{k,r}$-free.
Let $\lambda(G)$ be the spectral radius of $G$.
By the Perron--Frobenius Theorem,
we know that $\lambda(G)$
has an  eigenvector with all entries being positive,
we denote such an eigenvector by  $\mathbf{x}$.
For a vertex $v\in V(G)$,
we will write $\mathbf{x}_v$ for
the eigenvector entry of $\mathbf{x}$ corresponding to $v$.
We  may normalize $\mathbf{x}$ so that it has maximum entry equal to $1$, and let $z$ be a vertex such that $\mathbf{x}_z = 1$.
If there are multiple such vertices,
we choose and fix $z$ arbitrarily among them.

In the sequel,
we shall prove Theorem \ref{main result} iteratively, giving successively more precise estimates on both the structure of $G$ and the eigenvector entries of the vertices, until finally we can show that $e(G) = \mathrm{ex}(n, F_{k,r})$.

The proof of Theorem~\ref{main result}  is outlined as follows.

\begin{itemize}
\item[$\spadesuit$]
We apply Corollary \ref{coro22} to give a lower bound
 $e(G)\ge t_{r-1}(n) - o(n^2)$.
Moreover, $G$ has a very large multipartite subgraph
on parts $V_1,\ldots ,V_{r-1}$ such that $\frac{n}{r-1} - o(n)\le |V_i|
\le \frac{n}{r-1}+ o(n)$;
see Lemma \ref{approximate structure}.

\item[$\heartsuit$]
We show that
the number of vertices
that have $\Omega (n)$ neighbors in their own part is bounded by $o(n)$,  and the number of vertices
that have degree less than $(\frac{r-2}{r-1}- o(1))n$ is also bounded by
$o(n)$; see Lemmas \ref{lem size of W} and \ref{lem size of L} respectively.
Furthermore,
we will prove that
such vertices do not exist,
and each $G[V_i]$ is $K_{1,k}$-free and $M_{k}$-free;
see Lemmas \ref{WsubsetL} and \ref{eigenvector entries}.

\item[$\clubsuit$]
Based on the previous lemmas,
we shall refine the structure of $G$, and show that almost all vertices in $V_i$ are adjacent to every vertex in $V_i^c$, implying the presence of a large complete $(r-1)$-partite subgraph in $G$;
see Lemma \ref{Bi}.
Moreover, we shall prove that
$\mathbf{x}_u = 1- o(1)$ for every $u\in V(G)$;
see Lemma \ref{eigenvector}.

\item[$\diamondsuit$]
Once we know that all vertices have eigenvector entry  close to $1$, we can show that the $(r-1)$-partition is balanced;
see Lemma \ref{balance}.
Invoking these facts,
we finally show that $e(G) = \mathrm{ex}(n, F_{k,r})$.
\end{itemize}

\begin{lemma} \label{lem32}
Let $G$ be an $F_{k,r}$-free graph on $n$ vertices with maximum spectral radius. Then $$\lambda(G) \ge \bigg(1- \frac{1}{r-1}\bigg)n - \frac{r-1}{4n}.$$
\end{lemma}

\begin{proof}
 Let $H $
 be an $F_{k,r}$-free
 graph on $n$ vertices with maximum number of edges.
 Since $G$ is the graph maximizing the spectral radius over all $F_{k,r}$-free graphs, in view of Theorem \ref{thmChen},
 we can see by the Rayleigh quotient  that
\begin{equation*}
\lambda (G) \geq \lambda (H) \geq \frac{\mathbf{1}^T A(H) \mathbf{1}}{\mathbf{1}^T\mathbf{1}}
= \frac{ 2( t_{r-1}(n)
 + f(k-1,k-1))}{n}.
\end{equation*}
Note that $t_{r-1}(n) \ge (1-\frac{1}{r-1}) \frac{n^2}{2} - \frac{r-1}{8}$, so we have
$\lambda (G) \ge (1- \frac{1}{r-1})n - \frac{r-1}{4n}$.
\end{proof}

Applying Lemma \ref{lemniki} and Corollary \ref{coro22},
we obtain the  asymptotic structure of $G$.
Roughly speaking, we can find a large $(r-1)$-partite
subgraph in $G$.

\begin{lemma}[Approximate structure]
\label{approximate structure}
Let $G$ be an $F_{k,r}$-free graph on $n$ vertices with maximum spectral radius. For every $\epsilon > 0$, there is an integer $n_0$ such that if $n \geq n_0$, then $$e(G) \geq t_{r-1}(n) - \epsilon n^2.$$
Furthermore, there exists $\epsilon_1=\sqrt{6\epsilon}$ such that $G$ has a maximum $(r-1)$-cut $V = V_1 \cup \ldots \cup V_{r-1}$ with
\[ \sum_{1\leq i<j\leq r-1}e(V_i,V_j)\geq  t_{r-1}(n) - \epsilon n^2, \]
 and   for each $i \in [r-1]$,
 \[ \left(\frac{1}{r-1} - \epsilon_1\right)n \leq |V_i| \leq \left(\frac{1}{r-1} + \epsilon_1\right)n. \]
\end{lemma}

\begin{proof}
As suggested above, it follows from Lemma \ref{lemniki} and Corollary \ref{coro22} that for any given $\epsilon > 0$, we can take a large enough $n$ such that $e(G) \geq t_{r-1}(n) - \epsilon n^2.$  The same results also provide that there is a partition of $V(G) = U_1 \cup \ldots \cup U_{r-1}$ with $\sum_{i=1}^{r-1} e(U_i) \leq \epsilon n^2$, $\sum_{1\leq i<j\leq r-1}e(U_i,U_j)\geq  t_{r-1}(n) - \epsilon n^2$ and $\lfloor\frac{n}{r-1}\rfloor \leq |U_i| \leq \lceil\frac{n}{r-1}\rceil$ for each $i \in [r-1]$.
Thus, any maximum $(r-1)$-cut of $V = V_1 \cup \ldots \cup V_{r-1}$ must have $\sum_{i=1}^{r-1} e(V_i) \leq \sum_{i=1}^{r-1} e(U_i) \leq \epsilon n^2$ and $\sum_{1\leq i<j\leq r-1}e(V_i,V_j)\geq\sum_{1\leq i<j\leq r-1}e(U_i,U_j)\geq  t_{r-1}(n) - \epsilon n^2$.

Furthermore, since $G$ has edit distance at most $\epsilon n^2$ from some graph isomorphic to $T_{r-1}(n)$, we may let $a=\max\left\{\left||V_j|-\frac{n}{r-1}\right|, j\in [r-1]\right\}$.  Without loss of generality, we assume $\left||V_1|-\frac{n}{r-1}\right|=a$. Then
\begin{eqnarray*}
e(G)&\leq& \sum_{1\leq i<j\leq r-1}|V_i||V_j|+\sum_{i=1}^{r-1}e(V_i)\nonumber\\[2mm]
&=& |V_1|(n-|V_1|)+ \sum_{2\leq i<j\leq r-1}|V_i||V_j| +\epsilon n^2\nonumber\\
&=& |V_1|(n-|V_1|)+ \frac{1}{2}\Big((\sum_{j=2}^{r-1}|V_j|)^2-\sum_{j=2}^{r-1}|V_j|^2\Big) +\epsilon n^2\nonumber\\[2mm]
&\leq& |V_1|(n-|V_1|)+\frac{1}{2}(n-|V_1|)^2-\frac{1}{2(r-2)}(n-|V_1|)^2+\epsilon n^2\\[2mm]
&<&-\frac{r-1}{2(r-2)}a^2+\frac{r-2}{2(r-1)}n^2+\epsilon n^2,
\end{eqnarray*}
where  the last second inequality  holds by H\"{o}lder's inequality, and  the last inequality holds since $\left||V_1|-\frac{n}{r-1}\right|=a$. On the other hand, $$e(G)\geq t_{r-1}(n) - \epsilon n^2\geq (1-\frac{1}{r-1}) \frac{n^2}{2} - \frac{r-1}{8}-\epsilon n^2> \frac{r-2}{2(r-1)}n^2-2\epsilon n^2,$$
as $n$ is large enough. Therefore,
$\frac{r-1}{2(r-2)}a^2<3\epsilon n^2$, which implies that $a<\sqrt{\frac{6(r-2)\epsilon}{r-1}n^2}<\sqrt{6\epsilon}n= \epsilon_1 n$. The proof is completed.
\end{proof}

\begin{lemma}
\label{lem size of W}
Let $\epsilon$ and $\theta$ be two sufficiently small constants with $\epsilon < \theta^2/3$. We denote
\begin{equation}
\label{eqn defn of W}
    W := \cup_{i = 1}^{r-1} \{v \in V_i : |N_G(v) \cap V_i| \geq \theta n \}.
\end{equation}
For sufficiently large $n$, we have $$|W| \leq \frac{2 \theta}{3}n + \frac{2k^2}{\theta n} < \theta n.$$
\end{lemma}

\begin{proof}
We obtain from Lemma \ref{approximate structure} that  $\sum_{1\leq i<j\leq r-1}e(V_i,V_j)\geq  t_{r-1}(n) - \epsilon n^2.$ Hence,
\[
\sum_{i=1}^{r-1}e(V_i)=e(G)-\sum_{1\leq i<j\leq r-1}e(V_i,V_j)\leq t_{r-1}(n)+k^2-t_{r-1}(n) +\epsilon n^2\leq \epsilon n^2+k^2.
\]

On the other hand, if we let $W_i : = W \cap V_i$ for all $i \in [r-1]$, then
\[2 e(V_i) = \sum_{u\in V_i}d_{V_i}(u) \geq \sum_{u\in W_i}d_{V_i}(u) \geq |W_i|\theta n\]
Thus \[\sum_{i=1}^{r-1}e(V_i) \geq \sum_{i=1}^{r-1} \frac{|W_i|}{2}\theta n = \frac{|W|}{2}\theta n.\]
Therefore, we have that $\frac{|W|}{2}\theta n \leq \epsilon n^2 +k^2$. This proves that $|W| \leq \frac{2 \theta}{3}n +\frac{2k^2}{\theta n} < \theta n.$
\end{proof}

\begin{lemma}
\label{lem size of L}
Let $k \geq 2$ and $\frac{2(r-2)}{r-1}\epsilon <\epsilon_2^2
\ll \theta$. We denote
\begin{equation}
\label{eqn defn of L}
    L := \bigg\{v \in V(G) : d(v) \leq \bigg( 1 - \frac{1}{r-1} - \epsilon_2 \bigg)n\bigg\}.
\end{equation}
Then $|L| \leq \epsilon_3 n$, where $\epsilon_3 \ll \epsilon_2$
is a sufficiently small constant satisfying
$\frac{r-2}{2(r-1)} \epsilon_3^2 -\epsilon_2 \epsilon_3 + \epsilon <0$.
\end{lemma}

\begin{proof}
To prove this, assume to the contrary that the cardinality of $L$ is greater than $\epsilon_3 n$. Then there exists a subset $L' \subseteq L$ with $|L'| = \lfloor \epsilon_3 n \rfloor$.
Therefore,
\begin{equation*}
    \begin{split}
 e[G \setminus L'] \geq e(G) - \sum_{v \in L'}d(v)
 & \geq t_{r-1}(n) - \epsilon n^2 - \epsilon_3 n^2 \bigg(1 - \frac{1}{r-1} - \epsilon_2\bigg)\\
& > \frac{(n - \lfloor \epsilon_3 n \rfloor)^2}{2} \bigg(1-\frac{1}{r-1}\bigg) + k^2 \\
& \geq t_{r-1}(n - \lfloor \epsilon_3 n \rfloor) + k^2.
    \end{split}
\end{equation*}
However, this is a contradiction as the above lower bound for $e[G \setminus L']$ exceeds the upper bound on the number of edges in any $F_{k, r}$-free graph on $n - |L'|$ vertices.
\end{proof}

The following lemma was given in \cite{CFTZ20}.

\begin{lemma}[Cioab\u{a} et al. \cite{CFTZ20}] \label{inclusion exclusion lemma}
If $A_1, \ldots, A_p$ be finite sets, then
\[|A_1 \cap \ldots \cap A_p| \geq \sum_{i=1}^p |A_i| - (p-1)\bigg| \bigcup_{i=1}^p A_i \bigg|.\]
\end{lemma}

\begin{lemma} \label{WsubsetL}
Let $W$ and  $L$ be the sets of vertices defined in (\ref{eqn defn of W}) and (\ref{eqn defn of L}). Then $W \subseteq L$.
\end{lemma}

\begin{proof}
 Suppose on the contrary that there exists a vertex $u_0 \in W$ and $u_0 \notin L$.
   Without loss of generality, we may assume that $u_0\in V_1$.
   Since  $V_1,\ldots ,V_{r-1}$  form a maximum
   $(r-1)$-partite subgraph,
   we have $d_{V_1}(u_0)\le
d_{V_i}(u_0)$ for each $i\in [2,r-1]$. Indeed, otherwise,
   we can move the vertex $u_0$ into some part $V_i$ and strictly increase the number of edges between $V_1$ and $V_i$.
 Thus, we can get $d(u_0)\ge (r-1)d_{V_1}(u_0)$, which implies
 \[ d_{V_2}(u_0)\ge d(u_0) - d_{V_1}(u_0) - (r-3)n
 \left(\frac{1}{r-1} +
 \epsilon_1  \right).   \]
   On the other hand,
   invoking the fact that  $u_0 \not\in L$, we get $d(u_0)>
   (1-\frac{1}{r-1}-\epsilon_2 )n$. So
   \begin{align*}
    d_{V_2}(u_0) & \ge \left(1- \frac{1}{r-1}\right)
   d(u_0) - (r-3)n
 \left(\frac{1}{r-1} + \epsilon_1  \right)  \\
  & \ge \frac{n}{(r-1)^2} - \frac{r-2}{r-1} \epsilon_2 n
  - (r-3) \epsilon_1 n.
   \end{align*}
Recall from Lemmas \ref{lem size of W} and \ref{lem size of L} that
$  |W| <{\theta} n$ and $ |L|\le \epsilon_3n$ .
  Hence, for fixed $\theta, \epsilon_3$ and sufficiently large $n$, we have
\begin{equation*}
  |V_i\setminus (W\cup L)|\ge
  \left(\frac{1}{r-1}- \epsilon_1 \right)n-
 \theta n -
   \epsilon_3 n \ge  k.
\end{equation*}

{\bf Claim.} $u_0$ is adjacent to at most $k-1$ vertices in  $V_1\setminus (W\cup L)$.

\medskip

   Suppose that $u_0$ is adjacent to   $k$ vertices  $u_1^{(1)}, u_2^{(1)}, \ldots, u_{k}^{(1)}$ in $V_1\setminus (W\cup L)$.
   Since $u_j^{(1)}\not\in L$,  we have
   \[ d(u_j^{(1)})> \left(1-\frac{1}{r-1}-\epsilon_2 \right)n. \]
    On the other hand,
   we have $d_{V_1}(u_j^{(1)})< \theta n$ because
   $u_j^{(1)}\notin W$.
   So for each $j \in [k]$,
   \begin{align*}
    d_{V_2}(u_j^{(1)})&\ge d(u_j^{(1)})-d_{V_1}(u_j^{(1)})
   - (r-3) \left(\frac{1}{r-1}+ \epsilon_1 \right)n  \\
   & \ge \frac{n}{r-1} - \epsilon_2 n -\theta n - (r-3)\epsilon_1 n  .
  \end{align*}
   By Lemma~\ref{inclusion exclusion lemma},
   we consider the common neighbors of
   $u_0,u_1^{(1)}, \ldots ,u_k^{(1)}$ in $V_2$,
   \begin{eqnarray*}
&& \left| N_{V_2}(u_0) \cap N_{V_2}(u_1^{(1)})
\cap  \cdots \cap N_{V_2}(u_{k}^{(1)}) \setminus (W\cup L) \right| \\
 &\ge& d_{V_2}(u_0) +
    \sum_{j=1}^{k} d_{V_2}(u_j^{(1)})
    - k \left| V_2 \right| - |W| - |L| \\
 &\ge&  \tfrac{n}{(r-1)^2} - \tfrac{r-2}{r-1} \epsilon_2n
  - (r-3) \epsilon_1 n
  + k\left(\tfrac{n}{r-1} - \epsilon_2 n -\theta n - (r-3)\epsilon_1 n \right)
  - k\left(\tfrac{1}{r-1} + \epsilon_1 \right)n - \theta n - \epsilon_3 n \\
& \ge & \frac{n}{(r-1)^2} - o(n)   > k
\end{eqnarray*}
for sufficiently large $n$. So there exist $k$  vertices $u_1^{(2)},u_2^{(2)}, \ldots ,u_{k}^{(2)}$  in $V_2\setminus
(W\cup L)$ such that the   subgraph formed by two partitions $\{u_1^{(1)}, \ldots, u_{k}^{(1)}\}$ and $\{u_1^{(2)}, \ldots, u_{k}^{(2)}\}$ is a complete bipartite graph.
It is easy to see that the subgraph of $G$ formed by the vertex $u_0$
together with such a complete bipartite graph
can contain a copy of $F_{k,3}$ centered at the vertex $u_0$.
In the sequel, we shall extend this copy to the intersecting cliques
$F_{k,r}$.
Let $s\in [2,r-2]$ be a positive integer.
Assume that we have found the vertices
$u_1^{(i)},u_2^{(i)}, \ldots ,u_k^{(i)}\in V_i\setminus (W\cup L),
(i=1,2,\ldots ,s)$ such that
these vertices form a complete $s$-partite subgraph in $G$.
We next consider the common neighbors of these vertices
in $V_{s+1}$. Similarly, we get that for each $i\in [s]$ and $j\in [k]$,
   \begin{align*}
    d_{V_{s+1}}(u_j^{(i)})&\ge d(u_j^{(i)})-d_{V_i}(u_j^{(i)})
   - (r-3) \left(\frac{1}{r-1}+ \epsilon_1 \right)n  \\
   & \ge \frac{n}{r-1} - \epsilon_2 n -\theta n - (r-3)\epsilon_1 n  .
  \end{align*}
By Lemma~\ref{inclusion exclusion lemma} again, we can obtain
   \begin{eqnarray*}
&& \left| N_{V_{s+1}}(u_0) \cap \left(
\cap_{i\in [s],j\in [k]}N_{V_{s+1}}(u_j^{(i)}) \right)
\setminus (W\cup L) \right| \\
 &\ge& d_{V_{s+1}}(u_0) +
    \sum_{i\in [s], j\in [k]} d_{V_{s+1}}(u_j^{(i)})
    - ks \left| V_{s+1} \right| - |W| - |L| \\
 &\ge&  \tfrac{n}{(r-1)^2} - \tfrac{r-2}{r-1} \epsilon_2n
  - (r-3) \epsilon_1 n
  + ks\left(\tfrac{n}{r-1} - \epsilon_2 n -\theta n - (r-3)\epsilon_1 n \right)
  - ks\left(\tfrac{1}{r-1} + \epsilon_1 \right)n - \theta n - \epsilon_3 n \\
& \ge & \frac{n}{(r-1)^2} - o(n)   > k
\end{eqnarray*}
Thus we can find $k$ vertices
$u_1^{(s+1)},u_2^{(s+1)}, \ldots ,u_k^{(s+1)}\in V_{s+1}\setminus (W\cup L)$, which together with the previous vertices
$u_j^{(i)}\in V_i\setminus (W\cup L),
(i\in [s],j\in [k])$  form a complete $(s+1)$-partite subgraph in $G$.
Thus, for each $i \in [r-1]$, we can find $k$ vertices from every vertex part
 $V_i \setminus (W\cup L)$ such that these vertices together with
 $u_0$ form a copy of $F_{k,r}$ centered at $u_0$,
this is a contradiction.
Therefore $u_0$ is adjacent to at most $k-1$ vertices in  $V_1\setminus (W\cup L)$.

Hence, applying Lemmas \ref{lem size of W} and   \ref{lem size of L} again,
we have
   \begin{eqnarray*}
   d_{V_1}(u_0)&\le & |W|+|L|+k-1\\
   &<&   \frac{2\theta }{3} n +
 \frac{ 2k^2}{\theta n}  + \epsilon_3 n +k-1 \\
&<& \theta n
\end{eqnarray*}
for sufficiently large $n$.  This is a contradiction to the fact that $u_0\in W$.
Hence $W\subseteq L$.
 \end{proof}

\begin{lemma} \label{largeindependent}
For each $i$, there exists an independent set $I_i \subseteq V_i$ such that $$|I_i| \geq |V_i|-  \epsilon_3 n- k^2.$$
\end{lemma}

 \begin{proof}
 Since $V_i\setminus L$ is large enough
 by  Lemma \ref{lem size of L},
 we first prove that there exists a large complete multipartite subgraph
 between $V_1,V_2,\ldots ,V_{r-1}$.
 Let  $u_1^{(1)},u_2^{(1)}, \ldots, u_{2k}^{(1)} $ be $2k$ vertices
 chosen {\it arbitrarily} from $V_1 \setminus L$. Then $u_j^{(1)}\notin L$ which implies that
$d(u_{j}^{(1)})> \left(1- \frac{1}{r-1} - \epsilon_2 \right)n.$
Note that $W\subseteq L$ by Lemma~\ref{WsubsetL},
so $u_j^{(1)}\notin W$, then $d_{V_1}(u_j^{(1)})< \theta n$. Hence
 \begin{align*}
    d_{V_2}(u_j^{(1)})
   & \ge \frac{n}{r-1} - \epsilon_2 n -\theta n - (r-3)\epsilon_1 n  .
  \end{align*}
 Furthermore, by Lemma~\ref{inclusion exclusion lemma}, we  have
    \begin{eqnarray*}
&& \left| N_{V_2}(u_1^{(1)})  \cap N_{V_2}(u_2^{(1)})
\cap  \cdots \cap N_{V_2}(u_{2k}^{(1)}) \setminus L \right| \\
 &\ge&
    \sum_{j=1}^{2k} d_{V_2}(u_j^{(1)})
    - (2k-1) \left| V_2 \right|  - |L| \\
 &\ge&  2k\left(\tfrac{n}{r-1} - \epsilon_2 n -\theta n - (r-3)\epsilon_1 n \right)
  - (2k-1) \left(\tfrac{1}{r-1} + \epsilon_1 \right)n - \epsilon_3 n \\
& \ge & \frac{n}{r-1} - o(n)   > 2k
\end{eqnarray*}
 for  sufficiently  large $n$.
Hence there exist $2k$ vertices $u_1^{(2)},u_2^{(2)}, \ldots, u_{2k}^{(2)}
\in V_2$
such that the subgraph formed between the two parts
$\{u_1^{(1)}, \ldots, u_{2k}^{(1)}\}$ and
$\{u_1^{(2)},\ldots, u_{2k}^{(2)}\}$ is a complete bipartite graph.
Let $s\in [2,r-2]$ be a positive integer.
Assume that we have found the vertices
$u_1^{(i)},u_2^{(i)}, \ldots ,u_{2k}^{(i)}\in V_i\setminus L,
(i=1,2,\ldots ,s)$ such that
these vertices form a complete $s$-partite subgraph in $G$.
We next consider the common neighbors of these vertices
in $V_{s+1}$. Similarly, we get that for each $i\in [s]$ and $j\in [2k]$,
   \begin{align*}
    d_{V_{s+1}}(u_j^{(i)})
   & \ge \frac{n}{r-1} - \epsilon_2 n -\theta n - (r-3)\epsilon_1 n  .
  \end{align*}
By Lemma~\ref{inclusion exclusion lemma} again, we can obtain
   \begin{eqnarray*}
&& \left| \left(
\cap_{i\in [s],j\in [2k]}N_{V_{s+1}}(u_j^{(i)}) \right)
\setminus L \right| \\
 &\ge&   \sum_{i\in [s], j\in [2k]} d_{V_{s+1}}(u_j^{(i)})
    - (2ks-1) \left| V_{s+1} \right|   - |L| \\
 &\ge&   2ks \left(\tfrac{n}{r-1} - \epsilon_2 n  - (r-3)\epsilon_1 n \right)
  -  (2ks-1) \left(\tfrac{1}{r-1} + \epsilon_1 \right)n  - \epsilon_3 n \\
& \ge & \frac{n}{r-1} - o(n)   > 2k
\end{eqnarray*}
Thus we can find $2k$ vertices
$u_1^{(s+1)},u_2^{(s+1)}, \ldots ,u_{2k}^{(s+1)}\in V_{s+1}\setminus L$, which together with the  vertices
$u_j^{(i)}\in V_i\setminus L,
(i\in [s],j\in [2k])$  form a complete $(s+1)$-partite subgraph in $G$.
Thus, for any $2k$ vertices in $V_1\setminus L$,
we can find $2k$ vertices from
 $V_i \setminus L$ for each $i\in [2,r-1]$ such that all these vertices form a complete  $(r-1)$-partite subgraph in $G$.

\medskip
{\bf Claim.}  $G[V_1\setminus L]$ is both $K_{1, k}$-free
and $M_{k}$-free.
\medskip

Recall that  $G$ contains a large
complete $(r-1)$-partite subgraph with each part in $V_i\setminus L$.
If $G[V_1\setminus L]$ contains a copy of
$K_{1,k}$ centered at a vertex $u_0\in V_1$ with
leaves $u_1^{(1)},u_2^{(1)},\ldots ,u_{k}^{(1)}$, then
by the discussion above, we can embed the $F_{k,r}$ into $G$.
Therefore, $G[V_1\setminus L]$ is $K_{1,k}$-free.
Now,
we assume that
 $\{u_1^{(1)}u_2^{(1)},u_3^{(1)}u_4^{(1)}, \ldots ,
 u_{2k-1}^{(1)}u_{2k}^{(1)}\}$ is
  a matching of size $k$.
  Then for each $j\in [k]$,
  the vertices $u_{2j-1}^{(1)},u_{2j}^{(1)},u_1^{(2)},u_j^{(3)}\ldots ,u_j^{(r-1)}$
  form a clique of order $r$, and these $r$ cliques intersect at the vertex $u_1^{(2)}$.
  So $G[V_1\setminus L]$ is $M_{k}$-free.

Hence both the maximum degree and the maximum matching number of $G[V_1\setminus L]$  are at most  $k-1$, respectively. By Theorem \ref{Chvatal76},
$$e(G[V_1\setminus L])\le f(k-1, k-1).$$
 The same argument gives that for each $j\in [2,r-1]$,
 $$e(G[V_j\setminus L])\le f(k-1, k-1).$$
For each $i \in [r-1]$, since $G[V_i \setminus L]$ has at most $f(k-1, k-1)$ edges, then  the subgraph obtained from $G[V_i \setminus L]$ by deleting one vertex of each edge in $G[V_i\setminus L]$ contains no edges, which is an independent set of $G[V_i\setminus L]$.
 By Lemma \ref{lem size of L}, there exists an independent set
 $I_i\subseteq V_i$ such that
  \begin{align*}
  |I_i| &\ge |V_i\setminus L|-f(k-1, k-1) \ge |V_i|-  \epsilon_3 n- k^2 .
  \end{align*}
 This completes the proof.
 \end{proof}

\begin{lemma}
\label{eigenvector entries}
$L$ is empty, and each $G[V_i]$ is $K_{1,k}$-free and $M_k$-free.
\end{lemma}

\begin{proof}
 Recall that $A \mathbf{x} =\lambda (G) \mathbf{x}$ and
 $z$ is defined as a vertex with maximum eigenvector entry and satisfies $\mathbf{x}_z=1$. So we have
 $$d(z)\ge \sum_{w\sim z}\mathbf{x}_w=\lambda (G) \mathbf{x}_z=\lambda (G)\ge   \left(1-\frac{1}{r-1}- \frac{r-1}{4n^2} \right)n > \left(1-\frac{1}{r-1}- \epsilon_2 \right)n,$$
 as $n$ is large enough.
 Hence $z\notin L$.
  Without loss of generality, we may assume that $z\in V_1$.
Since the maximum degree in the induced subgraph $G[V_1\setminus L]$  is at most $k-1$ (containing no $K_{1,k}$),
from Lemma~\ref{lem size of L}, we have $|L|\le \epsilon_3 n$ and
 $$
 d_{V_1}(z)=d_{{V_1}\cap L}(z)+d_{V_1\setminus L}(z)\le
\epsilon_3 n + k-1 .
 $$
 Therefore,  by Lemma \ref{largeindependent}, we have
 \begin{eqnarray*}
 \lambda (G)
 &=&\lambda (G) \mathbf{x}_z=\sum_{v\sim z} \mathbf{x}_v
 = \sum_{\substack{v\in V_1 \\ v\sim z}} \mathbf{x}_v+
 \sum_{\substack{v\in V_2 \cup \cdots \cup V_{r-1} \\ v\sim z} } \mathbf{x}_v\\
 &=& \sum_{\substack{v\in V_1 \\ v\sim z}} \mathbf{x}_v+
 \sum_{\substack{v\in I_2 \cup \cdots \cup I_{r-1} \\ v\sim z}} \mathbf{x}_v+\sum_{\substack{v\in \cup_{i=2}^{r-1}V_i\setminus I_i \\ v\sim z} } \mathbf{x}_v\\
 &\le&  d_{V_1}(z)+\sum_{ v\in I_2\cup \cdots \cup I_{r-1}} \mathbf{x}_v+
 |\cup_{i=2}^{r-1}V_i\setminus I_i |\\
 &\le & \epsilon_3 n +k-1+\sum_{ v\in I_2\cup \cdots \cup I_{r-1}}
 \mathbf{x}_v+
 (r-2)(\epsilon_3n +k^2) .
\end{eqnarray*}
 By Lemma \ref{lem32}, we can get
 \begin{equation}\label{Lempty2}
 \sum_{ v\in I_2 \cup \cdots \cup I_{r-1}} \mathbf{x}_v\ge
\left(1-\frac{1}{r-1} - \frac{r-1}{4n^2}\right)n - (r-1)\epsilon_3 n - (r-2)k^2 - k+1.
 \end{equation}

Next we are going to prove   $L=\varnothing$.

By way of contradiction, assume that there is a vertex  $v\in L$,
so  $d_G(v)\le (1- \frac{1}{r-1} - \epsilon_2 )n$.
Consider the graph $G^+$ with vertex set $V(G)$ and edge set $E(G^+) = E(G \setminus \{v\}) \cup \{vw: w\in \cup_{i=2}^{r-1} I_i\}$.
Roughly speaking, in this process, the number of added edges is greater than
the number of deleted edges.
Note that adding a vertex incident with vertices in $I_i$ does not create any cliques,
and so $G^+$ is $F_{k,r}$-free.
Note that $\mathbf{x}$ is a vector such that
$\lambda (G)=\frac{\mathbf{x}^TA(G)\mathbf{x}}{
\mathbf{x}^T\mathbf{x}}$, and the Rayleigh theorem implies
$\lambda (G^+) \ge \frac{\mathbf{x}^TA(G^+)\mathbf{x}}{
\mathbf{x}^T\mathbf{x}}$.
Furthermore,
\begin{align*}
\lambda(G^+) - \lambda(G) &\geq \frac{\mathbf{x}^T\left(A(G^+) - A(G)\right) \mathbf{x}}{\mathbf{x}^T\mathbf{x}} = \frac{2\mathbf{x}_v}{\mathbf{x}^T\mathbf{x}}\left(
\sum_{w\in I_2\cup \cdots \cup I_{r-1}} \mathbf{x}_w - \sum_{uv\in E(G)} \mathbf{x}_u\right) \\
& \overset{(\ref{Lempty2})}{\geq}
\frac{2\mathbf{x}_v}{\mathbf{x}^T\mathbf{x}}
\left( \epsilon_2 n -\frac{r-1}{4n}-(r-1)\epsilon_3n - (r-2)k^2 -k+1 \right)>0,
\end{align*}
where the last inequality holds for $n$ large enough and
$\epsilon_3 \ll \epsilon_2$. This contradicts $G$ having the largest spectral radius over all $F_{k,r}$-free graphs, so $L$ must be empty.
Furthermore,
the claim in the proof of Lemma \ref{largeindependent}
implies that each $G[V_i]$ is $K_{1,k}$-free and $M_{k}$-free.
\end{proof}

\begin{lemma}\label{Bi}
For any $i\in [r-1]$, let $B_i=\{u\in V_i: d_{V_i}(u)\geq 1\}$ and $C_i=V_i\setminus B_i$. Then

(1) $|B_i|\leq 2k^2+1$;

(2) For every vertex $u\in C_i$, $u$ is adjacent to all vertices of $V\setminus V_i$.
\end{lemma}
\begin{proof}
We prove the  assertions by contradiction.

(1) If there exists a $j\in [r-1]$ such that $|B_j|> 2k^2+1$, then $\sum_{u\in B_j}d_{V_j}(u)>2k^2+1$. Since $G[V_j]$ is  both $K_{1,k}$-free and $M_k$-free, $e(G[V_j])\leq f(k-1,k-1)<k^2$. Therefore,
\[
2k^2+1< \sum_{u\in B_j}d_{V_j}(u)=\sum_{u\in V_j}d_{V_j}(u)=2e(G[V_j])<2k^2,
\]
which is a contradiction.

(2) If there exists a vertex $v\in C_1$ such that there is a vertex $w_{1,1}\notin V_1$ and $vw_{1,1}\notin E(G)$. Let $G'$ be the graph with $V(G')=V(G)$ and $E(G')=E(G)\cup \{vw_{1,1}\}$. We claim that $G'$ is $F_{k,r}$-free. Otherwise, $G'$ contains a copy of $F_{k,r}$, say  $F_0$, as a subgraph,
then $vw_{1,1}\in E(F_0)$. We may assume that $v$ is the center of $F_0$ (The case that $v$ is not the center of $F_0$ can be proved similarly). As $v$ is the center of $F_0$, there exist vertices $w_{1,1}, w_{1,2},\cdots, w_{1,r-1}, w_{2,1},\cdots, w_{2,r-1},\ldots,$ $w_{k,1},\cdots,w_{k,r-1}\notin V_1$ such that for any $i\in [k]$,
the vertex set $\{w_{i,1},w_{i,2},\ldots, w_{i,r-1}\}$ induces a copy of $K_{r-1}$ in $G$.  Therefore, for any $i\in [k]$ and $j\in [r-1]$, we have
\[
d_{V_1}(w_{i,j})=d(w_{i,j})-d_{V\setminus V_1}(w_{i,j})\geq \delta(G)-(k-1)-(r-3)\Big(\frac{n}{r-1}+ \epsilon_1 n\Big),
\]
where the last inequality holds as $G[V_s]$ is $K_{1,k}$-free, $|V_s|\leq \frac{n}{r-1}+\epsilon_1 n$ for any $s\in [r-1]$.
Since $L$ is empty by Lemma \ref{eigenvector entries}, we have $\delta(G)>(\frac{r-2}{r-1}-\epsilon_2)n$.
It follows  that
\[
d_{V_1}(w_{i,j})> \frac{n}{r-1}- o(n).
\]

\noindent
Using Lemma \ref{inclusion exclusion lemma}, we get
\begin{eqnarray*}
& & \Big|\bigcap_{i=1}^{k}\bigcap_{j=1}^{r-1}N_{V_1}(w_{i,j})\setminus B_1\Big|\\[2mm]
&\geq & \sum_{i=1}^{k}\sum_{j=1}^{r-1}|(N_{V_1}(w_{i,j})|-(k(r-1)-1)\Big|\bigcup_{i=1}^{k}\bigcup_{j=1}^{r-1}N_{V_1}(w_{i,j})\Big|- |B_1|\\[2mm]
&\geq &\sum_{i=1}^{k}\sum_{j=1}^{r-1}d_{V_1}(w_{i,j})-(kr-k-1)|V_1|- |B_1|\\[2mm]
&> & k(r-1)\left(\frac{n}{r-1}- o(n)\right)-
(kr-k-1)\left(\frac{n}{r-1}+o(n)\right)-(2k^2+1)\\[2mm]
&\geq & \frac{n}{r-1}-o(n) > 1.
\end{eqnarray*}
 Then there exists $v'\in C_1$ such that $v'$ is adjacent to $w_{1,1},\ldots,w_{1,r-1},\ldots, w_{k,1},\ldots,w_{k,r-1}$.
 Then $(F_0\setminus \{v\})\cup \{v'\}$ is a copy of  $F_{k,r}$ in $G$, which is a contradiction. Thus $G'$ is $F_{k,r}$-free. From the construction of $G'$, we see that $\lambda(G')>\lambda(G)$, which contradicts the assumption that $G$ has the maximum spectral radius among all $F_{k,r}$-free graphs on $n$ vertices.
\end{proof}

\begin{lemma}\label{eigenvector}
For any $u\in V(G)$,  $\mathbf{x}_u\geq 1-\frac{20k^2r^2}{n}$.
\end{lemma}

\begin{proof}
Recall that $\mathbf{x}_z=\max\{\mathbf{x}_i:i\in V(G)\}=1$. Without loss of generality, we may assume that $z\in V_1$. Then
\begin{align*}
\lambda (G) \mathbf{x}_z&=\sum_{w\thicksim z}\mathbf{x}_w = \sum_{w\thicksim z,w\in V_1}\mathbf{x}_w+\sum_{i=2}^{r-1}\Big(\sum_{w\thicksim z,w\in V_i}\mathbf{x}_w\Big)\\
&= \sum_{w\thicksim z,w\in V_1}\mathbf{x}_w+\sum_{i=2}^{r-1}\Big(\sum_{w\thicksim z,w\in B_i}\mathbf{x}_w+\sum_{w\thicksim z,w\in C_i}\mathbf{x}_w\Big),
\end{align*}
which implies that
\begin{align}
\sum_{i=2}^{r-1}\Big(\sum_{w\thicksim z,w\in C_i}\mathbf{x}_w\Big)
&=\lambda (G)-\sum_{w\thicksim z,w\in V_1}\mathbf{x}_w-\sum_{i=2}^{r-1}\Big(\sum_{w\thicksim z,w\in B_i}\mathbf{x}_w\Big)\nonumber\\
&\geq \lambda (G)-d_{V_1}(z)-\sum_{i=2}^{r-1}\Big(\sum_{w\in B_i}1\Big)\nonumber\\[2mm]
&\geq \lambda (G) -(k-1)-(r-2)(2k^2+1),\label{G'}
\end{align}
where  (\ref{G'}) holds as $G[V_1]$ is $K_{1,k}$-free, and $|B_i|\leq 2k^2+1$ for any $i\in [r-1]$.

We will prove this lemma by contradiction. Suppose that there is a vertex $v\in V(G)$ with $\mathbf{x}_v< 1-\frac{20k^2r^2}{n}$. Let $G'$ be the graph with $V(G')=V(G)$ and $E(G')=E(G\setminus \{v\})\cup \{vw: w\in N(z)\cap (\cup_{i=2}^{r-1}C_i)\}$. Since $C_i$ is an independent set for any $i\in [r-1]$, one may observe that $G'$ is $F_{k,r}$-free. By  (\ref{G'}), we have
\begin{align*}
\lambda(G')-\lambda(G)&\geq \frac{\mathbf{x}^{T}(A(G')-A(G))\mathbf{x}}{\mathbf{x}^T\mathbf{x}}\\[2mm]
&= \frac{2\mathbf{x}_v}{\mathbf{x}^T\mathbf{x}}\left(\sum_{i=2}^{r-1}\Bigl(\sum_{w\thicksim z,w\in C_i}\mathbf{x}_w \Bigr)-\sum_{uv\in E(G)}\mathbf{x}_u\right)\\[2mm]
&= \frac{2\mathbf{x}_v}{\mathbf{x}^T\mathbf{x}}\left(\sum_{i=2}^{r-1}\Bigl(\sum_{w\thicksim z,w\in C_i}\mathbf{x}_w \Bigr)-\lambda (G) \mathbf{x}_v\right)\\[2mm]
&>  \frac{2\mathbf{x}_v}{\mathbf{x}^T\mathbf{x}}\left(\lambda (G)-(k-1)-(r-2)(2k^2+1)-\lambda (G)\Bigl(1-\frac{20k^2r^2}{n}\Bigr)\right)\\[2mm]
&\geq  \frac{2\mathbf{x}_v}{\mathbf{x}^T\mathbf{x}}
\left(\frac{r-2}{r-1}20k^2r^2-\frac{r-1}{4n}\frac{20k^2r^2}{n}-k +1-(r-2)(2k^2+1)\right) >0,
\end{align*}
where the last  inequality follows by
$\lambda (G) \ge (1- \frac{1}{r-1})n - \frac{r-1}{4n}$ by Lemma \ref{lem32}.
 This contradicts the assumption that $G$ has the maximum spectral radius among all $F_{k,r}$-free graphs on $n$ vertices. Thus $\mathbf{x}_u\geq 1-\frac{20k^2r^2}{n}$ for any $u\in V(G)$.
\end{proof}

Let $G_{in}=\cup_{i=1}^{r-1}G[V_i]$. For any $i\in[r-1]$, let $|V_i|=n_i$ and $F=K_{n_1,n_2,\ldots,n_{r-1}}$ be the complete $(r-1)$-partite graph on $V_1,V_2,\ldots,V_{r-1}$. Let $G_{out}$ be the graph with $V(G_{out})=V(G)$ and $E(G_{out})=E(F)\setminus E(G)$.

\begin{lemma}\label{balance}
For any $1\leq i<j\leq r-1$,  $\left||V_i|-|V_j|\right|\leq 1$.
\end{lemma}

\begin{proof}
Suppose  $n_1\geq n_2\geq \ldots \geq n_{r-1}$. We prove the assertion  by contradiction.  Assume that there exist $i_0, j_0$ with $1\leq i_0 < j_0\leq r-1$ such that $n_{i_0}-n_{j_0}\geq 2$.

\noindent{\bfseries Claim 1.} There exists a constant $c_1>0$ such that $\lambda(T_{r-1}(n))-\lambda(F)\geq \frac{c_1}{n}$.

\begin{proof}
Let $F'=K_{n_1,\ldots, n_{i_0}-1,\ldots,n_{j_0}+1,\ldots,n_{r-1}}.$ Assume $F'\cong K_{n'_1,n'_2,\ldots,n'_{r-1}}$, where $n'_1\geq n'_2\geq \ldots \geq n'_{r-1}$.
By (\ref{eq0}), we have
\begin{equation}\label{C}
1=\sum_{i=1}^{r-1}\frac{n_i}{\lambda(F)+n_i}=\frac{n_{i_0}}{\lambda(F)+n_{i_0}}+\frac{n_{j_0}}{\lambda(F)+n_{j_0}}+\sum_{i\in [r-1]\setminus \{i_0,j_0\}}\frac{n_i}{\lambda(F)+n_i},
\end{equation}
and
\begin{equation}\label{C'}
1=\sum_{i=1}^{r-1}\frac{n'_i}{\lambda(F')+n'_i}=\frac{n_{i_0}-1}{\lambda(F')+n_{i_0}-1}+\frac{n_{j_0}+1}{\lambda(F')+n_{j_0}+1}+\sum_{i\in [r-1]\setminus\{i_0,j_0\}}\frac{n_i}{\lambda(F')+n_i}.
\end{equation}

Subtracting  (\ref{C'}) from   (\ref{C}), we get
\begin{eqnarray*}
& &\frac{2(n_{i_0}-n_{j_0}-1)\lambda^2(F)+(n_{i_0}+n_{j_0})(n_{i_0}-n_{j_0}-1)\lambda(F)}{(\lambda(F)+n_{i_0}-1)(\lambda(F)+n_{i_0})(\lambda(F)+n_{j_0}+1)(\lambda(F)+n_{j_0})}\\[2mm]
& = & \sum_{i\in [r-1]\setminus\{i_0,j_0\}}\frac{n_i(\lambda(F')-\lambda(F))}{(\lambda(F)+n_i)(\lambda(F')+n_i)}+\frac{(n_{i_0}-1)(\lambda(F')-\lambda(F))}{(\lambda(F)+n_{i_0}-1)(\lambda(F')+n_{i_0}-1)}\\[2mm]
& &+\frac{(n_{j_0}+1)(\lambda(F')-\lambda(F))}{(\lambda(F)+n_{j_0}+1)(\lambda(F')+n_{j_0}+1)}\\[2mm]
& \leq & \frac{\lambda(F')-\lambda(F)}{\lambda(F)+n'_{r-1}}\Big(\sum_{i\in [r-1]\setminus\{i_0,j_0\}}\frac{n_i}{\lambda(F')+n_i}+\frac{n_{i_0}-1}{\lambda(F')+n_{i_0}-1}+\frac{n_{j_0}+1}{\lambda(F')+n_{j_0}+1}\Big)\\[2mm]
& = & \frac{\lambda(F')-\lambda(F)}{\lambda(F)+n'_{r-1}},
\end{eqnarray*}
where the inequality holds as $n'_{r-1}\leq \min\{n_1,\ldots,n_{i_0}-1,\ldots,n_{j_0}+1,\ldots,n_{r-1}\}$, and the last equality is by  (\ref{C'}). Combining with the assumption $n_{i_0}-n_{j_0}\geq 2$, we obtain
\begin{eqnarray}\label{C-C'}
\frac{2\lambda^2(F)+(n_{i_0}+n_{j_0})\lambda(F)}{(\lambda(F)+n_{i_0}-1)(\lambda(F)+n_{i_0})(\lambda(F)+n_{j_0}+1)(\lambda(F)+n_{j_0})}\leq \frac{\lambda(F')-\lambda(F)}{\lambda(F)+n'_{r-1}}.
\end{eqnarray}
In view of  the construction of $F$, we see that
$$n-\Big(\frac{n}{r-1}+\epsilon_1 n\Big) \leq \delta(F)\leq \lambda(F)\leq \Delta(F)\leq n-\Big(\frac{n}{r-1}-\epsilon_1 n\Big),$$ thus $\lambda(F)=\Theta(n)$. From (\ref{C-C'}), it follows  that there exists a constant $c_1>0$ such that $\lambda(F')-\lambda(F)\geq \frac{c_1}{n}$. Therefore, by Lemma \ref{rpartitegraph}, $\lambda(T_{r-1}(n))-\lambda(F)\geq\lambda(F')-\lambda(F)\geq \frac{c_1}{n}$.

\end{proof}

\noindent{\bfseries Claim 2.} $$\lambda(G)\geq \lambda(T_{r-1}(n))+\frac{2f(k-1,k-1)}{n}\left(1-\frac{2}{n} \right).$$

\begin{proof}
Let $\mathbf{y}$ be an eigenvector of $T_{r-1}(n)$ corresponding to $\lambda(T_{r-1}(n))$,
$a=n-(r-1)\lfloor\frac{n}{r-1}\rfloor$.
Since $T_{r-1}(n)$ is a complete $(r-1)$-partite graph on $n$ vertices where
each partite set has either $\lfloor\frac{n}{r-1}\rfloor$ or $\lceil\frac{n}{r-1}\rceil$ vertices,  we may assume $\mathbf{y}=(\underbrace{\mathbf{y}_1,\ldots,\mathbf{y}_1}_{a\lceil\frac{n}{r-1}\rceil},\underbrace{\mathbf{y}_2,\ldots,\mathbf{y}_2}_{n-a\lceil\frac{n}{r-1}\rceil})^{\mathrm{T}}$.
Thus by  (\ref{eigenequation}), we have
\begin{eqnarray}
\lambda(T_{r-1}(n))\mathbf{y}_1=(r-a-1)\big\lfloor\frac{n}{r-1}\big\rfloor \mathbf{y}_2+(a-1)\big\lceil\frac{n}{r-1}\big\rceil \mathbf{y}_1,\label{14}
\end{eqnarray}
and \begin{eqnarray}
\lambda(T_{r-1}(n))\mathbf{y}_2=(r-a-2)\big\lfloor\frac{n}{r-1}\big\rfloor \mathbf{y}_2+a\big\lceil\frac{n}{r-1}\big\rceil \mathbf{y}_1.\label{15}
\end{eqnarray}
Combining (\ref{14}) and (\ref{15}), we obtain
\[
\Big(\lambda(T_{r-1}(n))+\big\lceil\frac{n}{r-1}\big\rceil\Big)
\mathbf{y}_1=\Big(\lambda(T_{r-1}(n))+\big\lfloor\frac{n}{r-1}\big\rfloor\Big)
\mathbf{y}_2.
\]
Without loss of generality, we assume that $\mathbf{y}_2=1$. Then
\[
\mathbf{y}_2\geq \mathbf{y}_1=\frac{\lambda(T_{r-1}(n))+\lfloor\frac{n}{r-1}\rfloor}{\lambda(T_{r-1}(n))+\lceil\frac{n}{r-1}\rceil}\geq 1-\frac{1}{\lambda(T_{r-1}(n))+\lceil\frac{n}{r-1}\rceil}.
\]
Since $\lambda(T_{r-1}(n))\geq \delta(T_{r-1}(n))\geq  n-\lceil\frac{n}{r-1}\rceil$,  $\mathbf{y}_1\geq 1-\frac{1}{n}$.

Let $H\in \mathrm{Ex}(n,F_{k,r})$. By Theorem \ref{thmChen}, $H$ is constructed from $T_{r-1}(n)$ by embedding a graph $G_0$ in one of the parts.
Then $e(H)=\mathrm{ex}(n,F_{k,r})=\mathrm{ex}(n,K_r)+f(k-1,k-1)$. Therefore
\begin{align}
\lambda(G)&\geq \lambda(H)\nonumber\geq \frac{\mathbf{y}^\mathrm{T}A(H)\mathbf{y}}{\mathbf{y}^{\mathrm{T}}\mathbf{y}}\nonumber\\[2mm]
&\geq \frac{\mathbf{y}^\mathrm{T}A(T_{r-1}(n))\mathbf{y}}{\mathbf{y}^{\mathrm{T}}\mathbf{y}}+\frac{2\sum_{ij\in E(G_0)}\mathbf{y}_{i}\mathbf{y}_{j}}{\mathbf{y}^{\mathrm{T}}\mathbf{y}}\nonumber\\[2mm]
&\geq \lambda(T_{r-1}(n))+\frac{2f(k-1,k-1)}{\mathbf{y}^{\mathrm{T}}\mathbf{y}}\left(1-\frac{1}{n}\right)^2\nonumber\\[2mm]
&\geq \lambda(T_{r-1}(n))+\frac{2f(k-1,k-1)}{n}\left(1-\frac{2}{n}\right).\label{Tnr-1}
\end{align}

\end{proof}

\noindent{\bfseries Claim 3.} $e(G_{in})-e(G_{out})\leq f(k-1,k-1).$

\begin{proof}
It follows from the definitions of  $G_{in}$ and $G_{out}$, that we have  $e(G_{in})=\sum_{i=1}^{r-1}|E(G[V_i])|$ and $e(G_{out})=\sum_{1\leq i<j\leq r-1}|V_i||V_j|-|E_{cr}(G)|$.
To get the claim, we need to prove (\ref{eqn1}) and (\ref{eqn2})  by Lemma \ref{chen}.
Obviously  (\ref{eqn2}) implies  (\ref{eqn1}), so it is sufficient to prove (\ref{eqn2}).
We  prove  (\ref{eqn2}) by contradiction. Without loss of generality, suppose that there exists a vertex $u\in V_1$ such that
$$d_{G[V_1]}(u)+\sum_{j=2}^{r-1}\beta(G[N(u)\cap V_j])\geq k.$$ Let $\{w_1w_2,\ldots,w_{2\ell-1}w_{2\ell}\}$ be an $\ell$-matching of $\cup_{j=2}^{r-1}G[N(u)\cap V_j]$ and $u_1,\ldots,u_{k-\ell}\in V_1$ be in the neighborhood of $u$. By Lemma \ref{Bi}, there exist $v_1,\ldots,v_{k-\ell}\in C_2$ such that $\{u,u_1,\ldots,u_{k-\ell},v_1,$ $\ldots,v_{k-\ell},w_1,\ldots,w_{2\ell}\}$ induce an $F_{k}$ of $G$.  For each $u_iv_i \ (1\leq i\leq k-\ell)$, there exist $r-3$ vertices $t_3\in C_3,t_4\in C_4,\ldots,t_{r-1}\in C_{r-1}$ such that $u,u_i,v_i,t_3,t_4,\ldots,t_{r-1}$ induce a $K_r$ of $G$. For any $ w_{i-1}w_i\in \{w_1w_2,\ldots,w_{2\ell-1}w_{2\ell}\}$, without loss of generality, suppose that $w_{i-1}w_i\subseteq E(G[V_2])$, then there exist $r-3$ vertices $z_3\in C_3,z_4\in C_4,\ldots,z_{r-1}\in C_{r-1}$ such that $u,w_{i-1},w_i,z_3,z_4,\ldots,z_{r-1}$ induce a $K_r$ of $G$. Thus we  find a copy of $F_{k,r}$ from the above $F_k$, a contradiction.

\end{proof}

According to the definitions of $G_{in}$, $G_{out}$ and $F$, we have $e(G)=e(G_{in})+e(F)-e(G_{out})$.
By Lemma \ref{Bi}, for any $i\in [r-1]$, and every vertex $u\in C_i$, $u$ is adjacent to all vertices of $V\setminus V_i$. Thus $$e(G_{out})\leq \sum_{1\leq i<j\leq r-1}|B_i||B_j|\leq \binom{r-1}{2}(2k^2+1)^2\leq 9k^4r^2.$$ Then
\begin{align}
\lambda(G)&= \frac{\mathbf{x}^\mathrm{T}A(G)\mathbf{x}}{\mathbf{x}^{\mathrm{T}}\mathbf{x}}\nonumber\\[2mm]
&= \frac{2\sum_{ij\in E(F)}\mathbf{x}_i\mathbf{x}_j}{\mathbf{x}^{\mathrm{T}}\mathbf{x}}+ \frac{2\sum_{ij\in E(G_{in})}\mathbf{x}_i\mathbf{x}_j}{\mathbf{x}^{\mathrm{T}}\mathbf{x}}- \frac{2\sum_{ij\in E(G_{out})}\mathbf{x}_i\mathbf{x}_j}{\mathbf{x}^{\mathrm{T}}\mathbf{x}} \nonumber\\[2mm]
&\leq \lambda(F)+ \frac{2e(G_{in})}{\mathbf{x}^{\mathrm{T}}\mathbf{x}}- \frac{2e(G_{out})(1-\frac{20k^2r^2}{n})^2}{\mathbf{x}^{\mathrm{T}}\mathbf{x}} \nonumber\\[2mm]
&\leq \lambda(F)+\frac{2(e(G_{in})-e(G_{out}))}{\mathbf{x}^{\mathrm{T}}\mathbf{x}}+\frac{2e(G_{out})\frac{40k^2r^2}{n}}{\mathbf{x}^{\mathrm{T}}\mathbf{x}}\nonumber\\[2mm]
&\leq \lambda(F)+\frac{2f(k-1,k-1)}{\mathbf{x}^{\mathrm{T}}\mathbf{x}}+\frac{\frac{720k^6r^4}{n}}{\mathbf{x}^{\mathrm{T}}\mathbf{x}}\label{GC}
\end{align}
Using   (\ref{Tnr-1}), (\ref{GC}) and $\mathbf{x}^{\mathrm{T}}\mathbf{x}\geq n(1-\frac{20k^2r^2}{n})^2\geq n-40k^2r^2$, we get
\begin{eqnarray*}
& &\lambda(T_{r-1}(n))-\lambda(F)\\[2mm]
&\leq & \frac{2f(k-1,k-1)}{\mathbf{x}^{\mathrm{T}}\mathbf{x}}-\frac{2f(k-1,k-1)}{n}+\frac{4f(k-1,k-1)}{n^2}+\frac{\frac{720k^6r^4}{n}}{\mathbf{x}^{\mathrm{T}}\mathbf{x}}\\[2mm]
&\leq & \frac{2f(k-1,k-1)}{n-40k^2r^2}-\frac{2f(k-1,k-1)}{n}+\frac{4f(k-1,k-1)}{n^2}+\frac{\frac{720k^6r^4}{n}}{n-40k^2r^2}\\[2mm]
&\leq & \frac{80k^2r^2f(k-1,k-1)}{n(n-40k^2r^2)}+\frac{4f(k-1,k-1)}{n^2}+\frac{720k^6r^4}{n(n-40k^2r^2)}\\[2mm]
&\leq & \frac{c_2}{n^2},
\end{eqnarray*}
where $c_2$ is a positive constant.

 Combining with Claim 1, we have
 \[
 \frac{c_1}{n}\leq \lambda(T_{r-1}(n))-\lambda(F)\leq \frac{c_2}{n^2},
 \]
 which is a contradiction when $n$ is sufficiently large. Thus $\left||V_i|-|V_j|\right|\leq 1$ for any $1\leq i<j\leq r-1$.

\end{proof}

\medskip
\noindent{\bfseries Proof of Theorem \ref{main result}.} Now we  prove that $e(G)=\mathrm{ex}(n,F_{k,r})$. Otherwise, we  assume that $e(G)\leq \mathrm{ex}(n,F_{k,r})-1$. Let $H$ be an $F_{k,r}$-free graph with $e(H)=\mathrm{ex}(n,F_{k,r})$ and $V(H)=V(G)$. By Lemma \ref{balance}, we may assume that $V_1,\ldots,V_{r-1}$ induce a complete $(r-1)$-partite graph in $H$. Let $E_1=E(G)\setminus E(H)$, $E_2=E(H)\setminus E(G)$,  then $E(H)=(E(G)\cup E_2)\setminus E_1$, and
\[
|E(G)\cap E(H)|+|E_1|=e(G)<e(H)=|E(G)\cap E(H)|+|E_2|,
\]
which implies that
$|E_2|\geq |E_1|+1$. Furthermore, by Lemma \ref{Bi}, we have
\begin{eqnarray}
|E_2|\leq f(k-1,k-1)+ \sum\limits_{1\leq i<j\leq r-1}|B_i||B_j|\leq k^2+\binom{r-1}{2}(2k^2+1)^2\leq 10k^4r^2.\label{eqn3}
\end{eqnarray}
According to   (\ref{Rayleigh}) and (\ref{eqn3}), we deduce, for sufficiently large $n$, that
\begin{align*}
\lambda(H)&\geq \frac{\mathbf{x}^{\mathrm{T}}A(H)\mathbf{x}}{\mathbf{x}^{\mathrm{T}}\mathbf{x}}\\[2mm]
&= \frac{\mathbf{x}^{\mathrm{T}}A(G)\mathbf{x}}{\mathbf{x}^{\mathrm{T}}\mathbf{x}}+\frac{2\sum_{ij\in E_2}\mathbf{x}_i\mathbf{x}_j}{\mathbf{x}^{\mathrm{T}}\mathbf{x}}-\frac{2\sum_{ij\in E_1}\mathbf{x}_i\mathbf{x}_j}{\mathbf{x}^{\mathrm{T}}\mathbf{x}}\\[2mm]
&= \lambda(G)+\frac{2}{\mathbf{x}^{\mathrm{T}}\mathbf{x}}\Big(\sum_{ij\in E_2}\mathbf{x}_i\mathbf{x}_j-\sum_{ij\in E_1}\mathbf{x}_i\mathbf{x}_j\Big)\\[2mm]
&\geq  \lambda(G)+\frac{2}{\mathbf{x}^{\mathrm{T}}\mathbf{x}}\Big(|E_2|(1-\frac{20k^2r^2}{n})^2- |E_1|\Big)\\[2mm]
&\geq  \lambda(G)+\frac{2}{\mathbf{x}^{\mathrm{T}}\mathbf{x}}\Big(|E_2|-\frac{40k^2r^2}{n}|E_2|- |E_1|\Big)\\[2mm]
&\geq \lambda(G)+\frac{2}{\mathbf{x}^{\mathrm{T}}\mathbf{x}}\Big(1-\frac{40k^2r^2}{n}|E_2|\Big)\\[2mm]
&\geq \lambda(G)+\frac{2}{\mathbf{x}^{\mathrm{T}}\mathbf{x}}\Big(1-\frac{40k^2r^2}{n}10k^4r^2\Big)\\[2mm]
&> \lambda(G),
\end{align*}
which contradicting the assumption that $G$ has the maximum spectral radius among all $F_{k,r}$-free graphs on $n$ vertices. Hence $e(G)=\mathrm{ex}(n,F_{k,r})$.
\qed

\section{Concluding remarks}

To avoid unnecessary calculations,
we did not attempt to get the best bound
on the order of graphs in the proof. 
 It would be  interesting
to determine how large $n$ needs to be for our result.

Recently,
Cioab\u{a}, Desai and Tait \cite{CDT21}
investigated the largest spectral radius of
 an $n$-vertex graph that does not contain the odd-wheel
 graph $W_{2k+1}$,
which is the graph obtained by joining a vertex to all vertices of
a cycle of length $2k$. Moreover, they raised
the following more general conjecture.

\begin{conjecture} \label{conj}
Let $F$ be any graph such that the graphs in $\mathrm{Ex}(n,F)$
are Tur\'{a}n graphs plus $O(1)$ edges.
Then for sufficiently large $n$,
a graph attaining the maximum spectral radius
among all $F$-free graphs on $n$ vertices is a member of $\mathrm{Ex}(n,F)$.
\end{conjecture}

We say that  $F$ is  edge-color-critical if
there exists an edge $e$ of $F$ such that
$\chi (F-e) < \chi (F)$.
Let $F$ be an edge-color-critical graph with $\chi (F)=r+1$.
By a result of Simonovits \cite{Sim66}
and a result of Nikiforov \cite{Niki09},
we know that $\mathrm{Ex}(n,F)=\mathrm{Ex}_{sp}(n,F) =
\{T_r(n)\}$ for sufficiently large $n$,
where $\mathrm{Ex}_{sp}(n,F)$ denotes the set of $F$-free graphs on $n$ vertices, attaining the maximum spectral radius.
This
 shows that Conjecture \ref{conj}  is true for
all edge-color-critical graphs.
As we mentioned before, Theorem \ref{thmCFTZ20}
says that Conjecture \ref{conj}
holds for the $k$-fan graph $F_k$.
Moreover, the result in \cite{Yongtao21} implies that
Conjecture \ref{conj} also holds
for the flower graph $H_{s,k}$, the graph defined by
intersecting $s$ triangles and $k$ odd cycles of length at least $5$
in exactly one common vertex.
In addition, our main result (Theorem \ref{main result})
tells us that    Conjecture \ref{conj} also holds
for the intersecting cliques $F_{k,r}$.
Note that $F_k, H_{s,k}$ and $F_{k,r}$ are not edge-color-critical.

Let  $S_{n,k}$ be the graph consisting of a clique on $k$ vertices and an independent set on $n-k$ vertices in which each vertex of the clique is adjacent to each vertex of the independent set.
Clearly, we can see that
$S_{n,k}$  does not contain $F_k$ as a subgraph.
Recently, Zhao, Huang and Guo \cite{ZHG21}
proved that $S_{n,k}$ is
the unique graph attaining the maximum signless Laplacian spectral radius among all graphs of order $n$ containing no $F_k$
for $n\ge 3k^2-k-2$.
Soon after, Chen, Liu and Zhang \cite{CLZ2021}
solved the corresponding case for $H_{s,k}$-free graphs.
So it is a natural question
to consider the maximum signless Laplacian spectral radius
among all graphs containing no $F_{k,r}$.
We write $q(G)$ for the
  signless Laplacian spectral radius, i.e.,
 the largest eigenvalue of
 the {\it signless Laplacian matrix} $Q(G)=D(G) +
 A(G)$, where $D(G)=\mathrm{diag}(d_1,\ldots ,d_n )$
 is the degree diagonal matrix and
 $A(G)$ is the adjacency matrix.
We end our paper with the following problem, and leave it for the interested readers.
Clearly, when $r=3$,
this problem reduces to the result of Zhao et al. \cite{ZHG21}.

\medskip
\noindent
{\bf Problem.}
{\it For integers $k\ge 1$ and $ r\ge 3$,
there exists an integer $n_0(k,r)$ such that
if $n\ge n_0(k,r)$ and
$G$ is an $F_{k,r}$-free graph on $n$
vertices, then
$  q(G) \le q(S_{n,k(r-2)})$,
equality holds if and only if $G=S_{n,k(r-2)}$.
}
\medskip

\section*{ Acknowledgements}
 We thank Erfang Shan and Yisai Xue for helpful suggestions.


\begin{thebibliography}{99}
\bibitem{Abbott72}
H.L. Abbott, D. Hanson,  H. Sauer, Intersection theorems for systems of sets,  J. Combin. Theory Ser.
A,  12 (1972) 381--389.




\bibitem{BG09}
L. Babai, B. Guiduli,
Spectral extrema for graphs: the Zarankiewicz problem,
Electronic J. Combin.,  15 (2009) R123.


\bibitem{Bapat14}
R.B. Bapat, Graphs and matrices, (2nd),
Universitext, Springer, London, 2014.

\bibitem{Bollobas78}
B. Bollob\'as, Extremal Graph Theory, Academic Press, New York, 1978.


\bibitem{Chen03}
G. Chen, R.J. Gould,
F. Pfender, B. Wei,
Extremal graphs for intersecting cliques,
 J. Combin. Theory  Ser. B,  89 (2003) 159--171.

  \bibitem{CLZ2021}
M.-Z. Chen, A.-M. Liu, X.-D. Zhang,
The signless Laplacian spectral radius of graphs without
intersecting odd cycles,
arXiv:2108.03895v1,  2021.

\bibitem{Chvatal76}
 V. Chv\'atal, D. Hanson, Degrees and matchings,  J. Combin. Theory Ser. B,  20  (1976) 128--138.

 \bibitem{CFTZ20}
 S. Cioab\u{a}, L.H. Feng, M. Tait, X.D. Zhang,
 The spectral radius of graphs with no intersecting triangles,
 Electron. J. Combin.,  27 (4) (2020) P4.22.


\bibitem{CDT21}
 S. Cioab\u{a}, D.N. Desai, M. Tait,
The spectral radius of graphs with no odd wheels,   arXiv: 2104.07729v1, 2021.

\bibitem{CDS1980}
D. Cvetkovi\'{c}, M. Doob, H. Sachs, Spectra of Graphs, Academic Press, New York, 1980.

\bibitem{Delorme}
C. Delorme, Eigenvalues of complete multipartite graphs, Discrete Math., 312 (2012)  2532-2535.

\bibitem{ES66}
P. Erd\H{o}os, M. Simonovits,
A limit theorem in graph theory, Stud. Sci. Math. Hungar.,  1 (1966) 51--57.

\bibitem{ES46}
P. Erd\H{o}s, A.H. Stone,
On the structure of linear graphs,
Bull. Am. Math. Soc., 52 (1946) 1087--1091.


\bibitem{Erdos95}
P. Erd{\H{o}}s, Z. F\"uredi, R.J. Gould, D.S. Gunderson,
Extremal Graphs for Intersecting Triangles,
 J. Combin. Theory  Ser. B,   64 (1995) 89--100.


 \bibitem{FLZ2007}
L.H. Feng, Q. Li, X.-D. Zhang,
Spectral radii of graphs with given chromatic number,
Appl. Math. Lett., 20 (2007)  158--162.


\bibitem{FiedlerNikif}
M. Fiedler, V. Nikiforov,  {Spectral radius and Hamiltonicity of graphs},   Linear Algebra
Appl.,  432 (2010)  2170--2173.






\bibitem{Furedi96}
Z. F\"{u}redi,
An upper bound on Zarankiewicz problem,
Comb. Probab. Comput., 5 (1996) 29--33.


\bibitem{Furedi96b}
Z. F\"{u}redi,
New asympotics for bipartite Tur\'an numbers,
J. Combin. Theory Ser. A, 75 (1996) 141--144.



\bibitem{Furedi2015} Z. F\"{u}redi, A proof of the stability of extremal graphs, Simonovits' stability from Szemer\'{e}di's regularity,
J. Combin. Theory Ser. B,   115 (2015) 66--71.

\bibitem{FS13}
Z. F\"uredi,  M. Simonovits,
The history of degenerate (bipartite) extremal graph problems,
in Erd\H{o}s Centennial,
Bolyai Soc. Math. Stud., 25,
J\'{a}nos Bolyai Math. Soc., Budapest, 2013, pp. 169--264.




\bibitem{Keevash11}
P. Keevash,
Hypergraph Tur\'{a}n problems, in Surveys in Combinatorics,
Cambridge University Press, Cambridge, 2011, pp. 83--140.


\bibitem{KST54}
T. K\"{o}v\'{a}ri, V.T. S\'{o}s, P. Tur\'{a}n,
On a problem of K. Zarankiewicz, Colloq. Math., 3 (1954) 50--57.


\bibitem{Man07}
W. Mantel, Problem 28, Solution by H. Gouwentak,
W. Mantel, J. Teixeira de
Mattes, F. Schuh and W. A. Wythoff.
Wiskundige Opgaven, 10 (1907) 60--61.

\bibitem{NikiforovTuran}
V. Nikiforov, Some inequalities for the largest eigenvalue of a graph,
 Combin. Probab. Comput.,  11  (2002)  179--189.

\bibitem{Nikiforov07}
V. Nikiforov,  Bounds on graph eigenvalues II,  Linear Algebra Appl., 427 (2007) 183--189.

\bibitem{Niki09}
V. Nikiforov,
Spectral saturation: inverting the spectral {T}ur\'{a}n theorem,
Electron. J. Combin., 16 (1) (2009) R33.

\bibitem{NikiforovKST}
V. Nikiforov, A contribution to the Zarankiewicz problem,
 Linear Algebra Appl., 414 (2010) 1405--1411.

\bibitem{NikiforovLAA10}
V. Nikiforov,  {The spectral radius of graphs without paths and cycles of specified length},  Linear Algebra Appl.,  432 (2010)  2243--2256.

\bibitem{Niki09b}
V. Nikiforov,
A spectral Erd\H{o}s-Stone-Bollob\'{a}s theorem,
Combin. Probab. Comput., 18  (3)  (2009)  455--458.

\bibitem{Niki09JGT}
V. Nikiforov, Stability for large forbidden subgraphs,
J. Graph Theory, 62   (4)  (2009)  362--368.

\bibitem{NikifSurvey}
V. Nikiforov,  Some new results in extremal graph theory,  Surveys in Combinatorics, London Math. Soc. Lecture Note Ser., 392, Cambridge Univ. Press, Cambridge, 2011, pp. 141--181.




\bibitem{Sim66}
M. Simonovits,
A method for solving extremal problems in graph theory, stability problems, in Theory of Graphs, Tihany, Hungary, 1966, Academic, New York, 1968, pp. 279--319.

\bibitem{Sim13}
M. Simonovits,
Paul Erd\H{o}s' influence on Extremal graph theory,
in The Mathematics of Paul Erd\H{o}s II, R.L. Graham, Springer, New York, 2013,
pp. 245--311.

\bibitem{Stevanovicetal} D. Stevanovi\'{c}, I. Gutman, M. Rehman, On spectral radius and enery of complete multipartite graphs, Ars Math. Contemp., 9 (2015) 109--113.

\bibitem{Turan41}
P. Tur\'{a}n,
On an extremal problem in graph theory,
Mat. Fiz. Lapok, 48 (1941)  436--452.
(in Hungarian).


\bibitem{Wilf86}
H.S. Wilf, Spectral bounds for the clique and independence numbers of graphs,  J. Combin. Theory Ser. B, 65  (1986) 113--117.



\bibitem{Yongtao21}
Y. Li, Y. Peng,
The spectral radius of graphs with no intersecting odd cycles,  arXiv:2106.00587v1, 2021.



\bibitem{ZW12}
M. Zhai, B. Wang,
Proof of a conjecture on the spectral radius of $C_4$-free
graphs, Linear Algebra Appl., 430 (2012) 1641--1647.



\bibitem{ZHG21}
Y. Zhao, X.Y. Huang, H. Guo,
The signless Laplacian spectral radius of graphs
with no intersecting triangles,
Linear Algebra Appl., 618 (2021) 12--21.

\end{thebibliography}
\end{document}